\documentclass[11pt,twoside,a4paper,leqno]{amsart}
\usepackage{amssymb}
\newif\ifdiags\diagstrue
\ifdiags\usepackage[PostScript=dvips,nohug]{diagrams}\fi
\theoremstyle{plain}
\newtheorem{thm}{Theorem}[section]
\theoremstyle{definition}
\newtheorem{defn}[thm]{Definition}
\newcommand{\emphdef}{\textit}
\numberwithin{equation}{section}

\newcommand{\thismonth}{\ifcase\month\or
  January\or February\or March\or April\or May\or June\or
  July\or August\or September\or October\or November\or December\fi
  \space\number\year}
\newcommand{\bauth}[1]{#1,}
\newcommand{\bart}[1]{\textit{#1},}
\newcommand{\bjourn}[3]{#1 #2 (#3)}
\newcommand{\bbook}[1]{#1,}
\newcommand{\bpp}[2]{#1--#2.}
\newcommand{\tcite}[1]{\textup{\cite{#1}}}
\DeclareMathAlphabet{\mathrmsl}{OT1}{cmr}{m}{sl}
\newcommand{\symb}[2]{\newcommand{#1}{\mathit{#2}} }
\newcommand{\rssymb}[2]{\newcommand{#1}{\mathrmsl{#2}} }
\symb\vol{vol}
\rssymb\iden{id}
\rssymb{\scal}{scal}
\newcommand{\R}{{\mathbb R}}
\newcommand{\C}{{\mathbb C}}
\newcommand{\eps}{\varepsilon}

\newcommand{\ip}[1]{\langle#1\rangle}

\newcommand{\tens}{\mathbin{\otimes}}
\newcommand{\dual}{^{*\!}}
\newcommand{\Cinf}{\mathrm{C}^\infty}
\newcommand{\SU}{\mathrm{SU}}
\newcommand{\connect}{\#}

\newcommand{\conf}{\mathsf{c}}
\newcommand{\cip}{\ip}
\newcommand{\cH}{{\mathcal H}}
\newcommand{\cS}{{\mathcal S}}
\newcommand{\CP}[1]{\C P^{#1}}
\newcommand{\gmw}{\mathrmsl w}
\newcommand{\gmv}{\mathrmsl v}
\begin{document}
\title[Einstein metrics and the Toda field equation]
{Einstein metrics, hypercomplex structures\\
and the Toda field equation}
\author{David M. J. Calderbank}
\address{Department of Mathematics and Statistics\\
University of Edinburgh\\ King's Buildings, Mayfield Road\\
Edinburgh EH9 3JZ\\ Scotland.}
\email{davidmjc@maths.ed.ac.uk}
\author{Paul Tod}
\address{Mathematics Institute\\ University of Oxford\\
24--29 St Giles\\ Oxford OX1 3LB\\ England.}
\email{tod@maths.ox.ac.uk}
\date{\thismonth}
\begin{abstract} We obtain explicitly all solutions of the
$\SU(\infty)$ Toda field equation with the property that the associated
Einstein-Weyl space admits a $2$-sphere of divergence-free shear-free geodesic
congruences. The solutions depend on an arbitrary holomorphic function and
give rise to new hyperK\"ahler and selfdual Einstein metrics with one
dimensional isometry group. These metrics each admit a compatible hypercomplex
structure with respect to which the symmetries are triholomorphic.
\end{abstract}
\maketitle

\section{Introduction}

LeBrun's $\cH$-space construction~\cite{LeBrun0} gives a method for
constructing locally, at least in principle, a selfdual Einstein manifold $M$
of negative scalar curvature with a prescribed conformal infinity given by an
arbitrary real analytic conformal $3$-manifold $B$. Furthermore, the geometry
of $M$ is uniquely determined by $B$. In practice, however, this ``filling
in'' construction is difficult to carry out directly, except in the case that
$B$ is conformally flat, when the Einstein metric is the hyperbolic metric. A
more general situation in which progress can be made is the case that $B$
admits a compatible Einstein-Weyl structure. Hitchin~\cite{Hitchin3} has shown
that in this situation, the twistor space $Z$ of $M$ is the projectivized
cotangent bundle of the minitwistor space $\cS$ of $B$, and consequently that
there is a conformal retraction of $M$ onto $B$, i.e., a conformal submersion
$M\to B$ inducing the identity map of $B$ at infinity. If one has enough
information about $B$ or its minitwistor space, then one can hope to find $M$
using this observation.

The first non-trivial examples of this construction were the Pedersen metrics
on the unit ball in $\R^4$~\cite{HP2}, where the conformal infinity is a
Berger sphere with its standard Einstein-Weyl structure~\cite{JT}. However,
even in these examples, $M$ is constructed indirectly and then shown to be the
desired $4$-metric using the uniqueness clause of the LeBrun construction. In
fact, Pedersen observes that the principal symmetry of the Berger sphere
induces a Killing field on $M$ and applies the Jones and Tod
construction~\cite{JT} to show that the Einstein metric on $M$ is given
locally by the total space of an abelian monopole over the round metric on
$S^3$. It is then not hard to guess which monopole does the job.

The following diagram summarizes the construction.
\ifdiags
\begin{diagram}[width=1.5em,height=1.5em,lefteqno,eqno=(1.0)]
    &                      &\;M^4&       &\\
    &\ldTo^\xi[crab=4pt]\ruIntoB&&\rdTo^K&\\
BS^3&                           &&       &S^3\\
\end{diagram}\fi
Here $K$ is a Killing field defining a Riemannian submersion, $\xi$ is a unit
vector field generating a conformal submersion, $S^3$ is the round $3$-sphere
and $BS^{3}$ is a Berger $3$-sphere. Note that we cannot take $M^4$ to be
the entire $4$-ball here, as $\xi$ is not globally defined on $B^4$ and $K$
vanishes at the origin. The conformal submersion $\xi$ induces an
Einstein-Weyl structure on $BS^3$. A Weyl structure may be specified by giving
a choice of representative $g$ for the conformal metric and a 1-form $\omega$
which defines a connection on the line bundle $L^1$ of scalars of weight one
(see below). On the Berger spheres, $g$ and $\omega$ are explicitly
\begin{align*}
g &= d\theta^2+\sin^2\theta d\phi^2+a^2(d\psi+\cos\theta d\phi)^2
=\sigma_1^{\,2}+\sigma_2^{\,2}+a^2\sigma_3^{\,2}\\
\omega &= b(d\psi+\cos\theta d\phi) = b\sigma_3.
\end{align*}
where $a,b$ are constants with $b^2=a^2(1-a^2)$ and the $\sigma_i$ are the
usual invariant 1-forms on $\SU(2)$. The Killing field $K$ on $M^4$ is induced
by the principal symmetry $\partial/\partial\psi$ on $BS^3$ and the Einstein
metric on $M^4$ is given by
\begin{equation*}\notag
g_M=\frac1{(1-\rho^2)^2}\left[\frac{1+m^2\rho^2}{1+m^2\rho^4}d\rho^2
+\frac14\rho^2\left((1+m^2\rho^2)(\sigma_1^{\,2}+\sigma_2^{\,2})
+\frac{1+m^2\rho^4}{1+m^2\rho^2}\sigma_3^2\right)\right].
\end{equation*}
This metric is complete on the unit ball, and the conformal structure extends
to the conformal infinity at $\rho=1$, which is the Berger sphere with
$a^2=1/(1+m^2)$.

Our aim in this paper is to generalize these metrics and the diagram (1.0) by
replacing $S^3$ and $BS^{3}$ with other Einstein-Weyl spaces. To do this, we
want to explain how the geometry of $M^4$ restricts the possible geometries of
$M^4/K$. First of all, by the general theory of~\cite{Tod3}, since $M^4$ is an
Einstein manifold with Killing field $K$, the Einstein metric is conformal to
a scalar-flat K\"ahler metric. Now an Einstein-Weyl space arising as the
quotient of a scalar-flat K\"ahler metric by a Killing field is not
arbitrary~\cite{LeBrun1}: it admits a shear-free twist-free geodesic
congruence~\cite{Tod2}. Let us pause to define some of these terms.
\begin{defn} A \emphdef{Weyl space} is a conformal manifold $(B^n,\conf)$
equipped with a torsion free connection $D$ with $D\conf=0$. (Here we view the
conformal structure $\conf$ as a metric on $TB$ with values in the real line
bundle $L^2$, where $L^{-n}=|\Lambda^nT^*B|$.) It is said to be
\emphdef{Einstein-Weyl} iff the symmetric tracefree part of the Ricci tensor
of this connection vanishes.

A \emphdef{congruence} on an oriented three-dimensional Weyl space $B^3$ is
generated by a weightless unit vector field $\chi\in \Cinf(B^3,L^{-1}\tens
TB)$, i.e., $\cip{\chi,\chi}=1$, where the angle brackets denote the conformal
metric. The congruence is \emph{shear-free and geodesic} iff
$$D\chi=\tau(\iden-\chi\tens\chi)+\kappa\,{*\chi}$$
and $\tau,\kappa$ are called the \emphdef{divergence} and \emphdef{twist}
of $\chi$. They are sections of $L^{-1}$.
\end{defn}
\begin{defn}
We shall say an Einstein-Weyl $3$-manifold is \emphdef{Toda} iff it admits a
shear-free geodesic congruence with vanishing twist, and that it is
\emphdef{hyperCR} iff it admits a shear-free geodesic congruence with
vanishing divergence.
\end{defn}
The reason for this terminology is that if $B^3$ is Toda, there is a
distinguished metric, which we call the \emphdef{LeBrun-Ward gauge},
such that the Weyl structure may be written
\begin{equation}\label{TodaForm}\begin{split}
g &= e^u(dx^2+dy^2)+dz^2\\ \omega &= -u_zdz,
\end{split}\end{equation}
where $u$ is a solution of the Toda field equation
$u_{xx}+u_{yy}+(e^u)_{zz}=0$ (see~\cite{Ward}).

As before, $\omega$ is the
connection $1$-form of the covariant derivative on $L^1$ given by the Weyl
connection, relative to the trivialization of $L^1$ given by the choice of
representative metric. Hence $Dg=-2\omega\tens g$.

On the other hand if $B^3$ is hyperCR, then $\chi$ is not alone: in fact
$D^B-\kappa\,{*1}$ is a flat connection on $L^{-1}\tens TB$ such that the
parallel weightless unit vector fields give a $2$-sphere of divergence-free
shear-free geodesic congruences. Each of these congruences defines a CR
structure on $B^3$, hence the term hyperCR is introduced by analogy with
hypercomplex or hyperK\"ahler. These Einstein-Weyl spaces were called
\emphdef{special} in~\cite{CTV} and~\cite{GT}. It was shown there that
monopoles over hyperCR Einstein-Weyl spaces define hypercomplex 4-manifolds by
the construction of~\cite{JT}.

The round metric on $S^3$ is Einstein, and therefore Einstein-Weyl. As such,
it admits a Toda congruence on the complement of any pair of antipodal points,
and also two families of hyperCR congruences. The Toda congruence is given by
the geodesics joining the antipodal points, while the two hyperCR structures
are given by the left and right invariant congruences respectively.

Returning now to the Pedersen metric over $S^3$, we see that not only is it
conformally scalar-flat K\"ahler, but also, since $S^3$ is hyperCR, it admits
two compatible hypercomplex structures~\cite{GT,Madsen2}. We shall see later
that these hypercomplex structures are also induced by a hyperCR structure on
the Berger $3$-sphere at infinity.

Therefore, in order to generalize the Pedersen metrics and diagram (1.0), one
approach is to look for Einstein-Weyl spaces, generalizing $S^3$, which are
both Toda and hyperCR. Our first result is that all such spaces can be found:

\begin{thm}\label{Main} For any holomorphic function $h$
on an open subset of $S^2$, the Einstein-Weyl space given by
\begin{align*}
g&=(z+h)(z+\overline h)g_{S^2}+dz^2\\
\omega&=-\frac{2z+h+\overline h}{(z+h)(z+\overline h)}dz,
\end{align*}
where $g_{S^2}$ is the spherical metric, is both Toda and hyperCR.
Furthermore any Toda Einstein-Weyl space admitting a hyperCR structure arises
in this way, with the exception of the Toda solutions given by a parallel
congruence on flat space.
\end{thm}

We prove this theorem in section~\ref{solns}. Then, in sections~\ref{sfk}
and~\ref{ein}, we consider monopoles over these hyperCR Toda spaces. In
particular, over each such space, we find an Einstein metric with symmetry and
with a conformal infinity given by another Einstein-Weyl space from a class
generalizing the Berger spheres.  This will give the desired generalization of
diagram (1.0).

\subsection*{Acknowledgements}

We would like to thank Ian Strachan for drawing our attention to the
linearized solutions of the Toda equation given in section~\ref{sfk}.
Part of this work was carried out while the first author was visiting
Humboldt-Universit\"at zu Berlin and the second author was visiting
the University of Adelaide. The authors gratefully acknowledge
support from SFB 288 and the Australian Research Council.
The diagrams were typeset with Paul Taylor's commutative diagrams
package.

\section{The Toda solutions}\label{solns}

In order to prove Theorem~\ref{Main} we must find all solutions of the Toda
field equation admitting a hyperCR structure. As very few solutions of the
Toda equation are known, this is an interesting exercise in its own right.
The condition that an Einstein-Weyl space is hyperCR is equivalent~\cite{GT}
to the existence of a section $\kappa$ of $L^{-1}$ with
\begin{align}
\kappa^2&=\frac16\scal^D\label{scaleqn}\\
D\kappa&=-\frac12{*F^D},\label{fareqn}
\end{align}
where $F^D$ is the curvature of the Weyl connection on $L^1$.  Our aim
is to impose this condition on the Toda field equation.  We start with
equation~\eqref{fareqn} which can be written in a gauge as
$d\kappa-\omega\kappa=-\frac12{*d\omega}$. In the LeBrun-Ward gauge, this
becomes
$$\kappa_xdx+\kappa_ydy+(\kappa_z+u_z\kappa)dz
=-\frac12u_{yz}dx+\frac12u_{xz}dy.$$ From this we deduce the equations
$\kappa_x=-\frac12u_{yz}$ and $\kappa_y=\frac12 u_{xz}$ which have an
integrability condition:
$$u_{xxz}=u_{xzx}=2\kappa_{yx}=2\kappa_{xy}=-u_{yzy}=-u_{yyz}.$$
Therefore $0=(u_{xx}+u_{yy})_z=-(e^u)_{zzz}$ and so we may write
$$e^u=e^{f(x,y)}\bigl(az^2+b(x,y)z+c(x,y)\bigr).$$
The Toda field equation with this Ansatz can be solved explicitly as
follows. We compute:
\begin{align*}
u_{xx}+u_{yy}&=\frac{\bigl((b_{xx}+b_{yy})z+c_{xx}+c_{yy}\bigr)(az^2+bz+c)
-(b_xz+c_x)^2-(b_yz+c_y)^2}{(az^2+bz+c)^2}\\
&\quad+f_{xx}+f_{yy}\\
(e^u)_{zz}&=2ae^f,
\end{align*}
which must sum to zero. We multiply through by $(az^2+bz+c)^2$ and equate
coefficients of the resulting quartic in $z$. The leading term is a Liouville
equation:
$$f_{xx}+f_{yy}+2ae^f=0.$$
The general solution of this Liouville equation may be written
$$e^{f(x,y)}=\frac{4|F'(x+iy)|^2}{(1+a|F(x+iy)|^2)^2}$$
in terms of an arbitrary nonconstant $F$, holomorphic in $x+iy$. The other
coefficients now give the following equations:
\begin{align}
a(b_{xx}+b_{yy})&=0\label{one}\\
a(c_{xx}+c_{yy})+b(b_{xx}+b_{yy})&=b_x^2+b_y^2\label{two}\\
b(c_{xx}+c_{yy})+c(b_{xx}+b_{yy})&=2(b_xc_x+b_yc_y)\label{three}\\
c(c_{xx}+c_{yy})&=c_x^2+c_y^2.\label{four}
\end{align}

If $a=0$ then equations~\eqref{two} and~\eqref{four} are solved by letting
$b=B|e^\phi|^2$, $c=C|e^\psi|^2$ with $B,C$ constant and $\phi,\psi$
holomorphic. Equation~\eqref{three} now gives $B=0$, $C=0$ or
$|\phi'-\psi'|=0$ and so the functional dependence of $b$ and $c$ can be
absorbed into $f$ and we have a separable solution for $e^u$.

If $a$ is not zero, then equations~\eqref{one} and~\eqref{four} give
$b=a(h+\overline h)$, $c=C|e^\psi|$, with $C$ constant and $h,\psi$
holomorphic. Equation~\eqref{two} now gives $aC|\psi'|^2|e^\psi|^2=a^2|h'|^2$
and so $C/a$ is nonnegative. If $C=0$ then $h$ is constant and we have a
separable solution; otherwise, without loss of generality, we may take $C=a$
and $h=e^\psi+\mu$ where $\mu$ is a real constant. Finally,
equation~\eqref{three} reduces to $\mu|\psi'|^2=0$ and so either $\mu=0$ or
$\psi$ is constant, the latter case again giving a separable solution.

The separable solutions are all known and the Einstein-Weyl structures are
all given by $3$-metrics of constant curvature~\cite{Tod3}: in our case
the curvature must be nonnegative in order to satisfy~\eqref{scaleqn},
and these solutions, generating the metrics of $\R^3$ and $S^3$, are the
ones we are trying to generalize. The new solutions of the Toda equation
are:
$$e^u=\frac{4a(z+h)(z+\overline h)|F'|^2}{(1+a|F|^2)^2},$$ and positivity
forces $a>0$. We readily verify that the remaining equations for $\kappa$
are satisfied with
$$\kappa=\frac{i(h-\overline h)}{2(z+h)(z+\overline h)},$$
and so the Einstein-Weyl space is indeed hyperCR.  Since $F$ cannot be
constant, we may use $\sqrt a F$ as a holomorphic coordinate in place of
$x+iy$ and Theorem~\ref{Main} easily follows.

\section{Scalar-flat K\"ahler metrics}\label{sfk}

In this section we study abelian monopoles over the hyperCR Toda spaces. On
any Toda space $B$, the monopole equation in the LeBrun-Ward gauge is
$$\gmw_{xx}+\gmw_{yy}+(e^u\gmw)_{zz}=0.$$ LeBrun~\cite{LeBrun1} shows that
each solution of this equation generates a scalar-flat K\"ahler metric with a
Killing field, given explicitly by
$$g_M = \gmw\,e^u(dx^2+dy^2)+ \gmw dz^2+ \gmw^{-1}(dt+\theta)^2,$$
where $\theta$ is a $1$-form on $B$ with $*(d\gmw-\omega\gmw)=d\theta$.

Consequently we can construct a large family of scalar-flat K\"ahler metrics
from the Einstein-Weyl spaces of section~\ref{solns}. Since the Einstein-Weyl
spaces are hyperCR, these scalar-flat K\"ahler metrics admit compatible
hypercomplex structures such that the symmetry $\partial/\partial t$ is
triholomorphic. For most choices of $h$, the hyperCR Toda spaces have no
continuous symmetries, and so these scalar-flat K\"ahler spaces will
generically have only a one-dimensional symmetry group, generated by
$\partial/\partial t$.

In order to obtain explicit metrics, we still have a linear differential
equation, the monopole equation, to solve. Fortunately, there are some
interesting solutions available, given to us for free by the geometry. These
solutions may be viewed as arising from LeBrun's observation~\cite{LeBrun1}
that the monopole equation above is the linearized Toda equation, and so
monopoles can be found by linearizing a family of solutions of the Toda
equation. In particular, the affine change $(x,y,z)\mapsto (ax,ay,az-b)$
induces a symmetry of the Toda equation.  Linearizing a family of solutions
generated by this gauge transformation shows that for any $a,b\in\R$,
$a(1-\frac12zu_z)+\frac12bu_z$ defines a monopole on \emphdef{any} Toda
Einstein-Weyl space~\cite{BF,GD,LeBrun1,Tod3}.

For our explicit solutions, these monopoles may also be obtained by
linearizing with respect to affine changes of the holomorphic function $h$.
Ian Strachan (private communication) has pointed out that by linearizing the
solutions with respect to arbitrary holomorphic changes of $h$, one sees, more
generally, that
$$\gmw=\frac{f}{2(z+h)}+\frac{\overline f}{2(z+\overline h)}$$ is a monopole
for any holomorphic function $f$ ($f=ah+b$ being a special case). To compute
$\theta$ note that $*(d\gmw-\omega\gmw)=d(\gmv\,dz)+
\frac12(f+\overline f)\vol_{S^2}$ where
$$\gmv=\frac{if}{2(z+h)}-\frac{i\overline f}{2(z+\overline h)}.$$
Hence one can write $dt+\theta=\beta+\gmv\,dz$, where $\beta$ is a
$1$-form independent of $z$ such that
$d\beta=\frac12(f+\overline f)\vol_{S^2}$,
so that the scalar-flat K\"ahler metric is:
\begin{equation}\label{sfkm}
g_M = \gmw\,(z+h)(z+\overline h)g_{S^2} + \gmw dz^2 +
\gmw^{-1}(\beta+\gmv\,dz)^2.
\end{equation}
For definiteness, one could take
$$\beta=dt+\frac i{1+\zeta\overline\zeta}\left(\frac{f\,d\zeta}{\zeta}-
\frac{\overline f\,d\overline\zeta}{\overline\zeta}\right)$$ where $\zeta$ is
a holomorphic coordinate on $S^2$ with $\vol_{S^2}=2i\,d\zeta\wedge
d\overline\zeta/(1+\zeta\overline\zeta)^2$.

These scalar-flat K\"ahler metrics will not be Einstein or conformally Einstein
in general. However, they do have the property that the lift of
$\partial/\partial z$ given by $\beta(\partial/\partial z)=0$ defines a
conformal submersion.  To see this, write the conformal structure on $M$ as
$\conf=(\eps_0)^2+\cdots+(\eps_3)^2$, where $\eps_0$ and $\eps_3$ are the
weightless unit 1-forms corresponding to $\gmw dz$ and $\beta+\gmv\,dz$. Let
$\xi$ be the weightless unit $1$-form dual to $\partial/\partial z$, so that
$$\xi=\frac{\gmw\eps_0+\gmv\eps_3}{\sqrt{\gmw^2+\gmv^2}}.$$  Now
$\eps_0^2+\eps_3^2-\xi^2=\eta^2$, where
$$\eta=\frac{\gmv\eps_0-\gmw\eps_3}{\sqrt{\gmw^2+\gmv^2}}
=\frac{\gmw\beta}{\sqrt{\gmw^2+\gmv^2}}.$$ Hence $\conf-\xi^2$ may be
represented by the metric
$$(\gmw^2+\gmv^2)|z+h|^2g_{S^2}+\beta^2=|f|^2g_{S^2}+\beta^2,$$
which is independent of $z$, so that $\xi$ is a conformal submersion
over this metric.

The conformal structures $|f|^2g_{S^2}+\beta^2$, depending on an
arbitrary holomorphic function $f$, arise elsewhere, namely in the
classification of Einstein-Weyl spaces admitting a continuous symmetry
(a conformal vector field preserving the Weyl connection) with
geodesic orbits:
\begin{thm}\tcite{CP}\label{gs}
The three dimensional Einstein-Weyl spaces with geodesic symmetry are either
flat with translational symmetry or are given locally by:
\begin{align*}
g&=|H|^{-2}g_{S^2}+\beta^2\\
\omega&=\tfrac i2(H-\overline H)\beta\\
d\beta&=\tfrac12(H+\overline H)|H|^{-2}\vol_{S^2}
\end{align*}
where $H$ is any nonvanishing holomorphic function on an open subset of $S^2$.
The geodesic symmetry $K$ is dual to $\beta$ and the monopoles associated to
$K/|K|$ are $\tau=\frac i2(H-\overline H)\mu_g^{-1}$,
$\kappa=\frac14(H+\overline H)\mu_g^{-1}$. Furthermore, these spaces are all
hyperCR, the flat connection on $L^{-1}\tens TB$ being $D^B+\kappa\,{*1}$.
\end{thm}

Replacing $H$ by $1/f$ we see that the scalar flat K\"ahler metrics of this
section fibre over the Einstein-Weyl spaces with geodesic symmetry. In
the next section we shall fill in these Einstein-Weyl spaces with Einstein
metrics.

\section{Selfdual Einstein metrics}\label{ein}

In the previous section we noted in passing, that when $f=ah+b$, the monopoles
$\gmw=\frac12\bigl(f/(z+h)+\overline f/(z+\overline h)\bigr)$ may be identified
with the geometrically significant monopoles $a(1-\frac12zu_z)+\frac12bu_z$
which are canonically defined on any Toda Einstein-Weyl space. There is also a
special monopole defined on any hyperCR space~\cite{CTV,GT}, namely $\kappa$,
as one easily sees from equation~\eqref{fareqn}. We note that for the hyperCR
Toda spaces, this monopole is obtained by setting $f=i$.

The significance of these monopoles is that they all give rise to Einstein
metrics.

\subsection{Scalar-flat K\"ahler metrics which are conformally hyperK\"ahler}
\noindent\newline
By~\cite{CTV,GT}, the hypercomplex structure we obtain from the $\kappa$
monopole is conformally hyperK\"ahler, and the symmetry $\partial/\partial t$
is a triholomorphic homothetic vector field of the hyperK\"ahler metric.
Hence when $f=i$, the scalar-flat K\"ahler metric~\eqref{sfkm} is
conformally Ricci-flat. The Einstein-Weyl space with geodesic
symmetry obtained from the conformal submersion $\xi$ is $\R^3$ (with
radial symmetry).

\subsection{HyperK\"ahler metrics with compatible hypercomplex structures}
\label{HK}\noindent\newline
By~\cite{BF,GD,LeBrun1}, on any Toda space, the scalar-flat K\"ahler metric
corresponding to the monopole $u_z$ is in fact Ricci-flat and therefore
hyperK\"ahler: the symmetry $\partial/\partial t$ is a Killing field of the
hyperK\"ahler metric, but is not triholomorphic unless the Toda space is
$\R^3$ (with a translational congruence). However, in the case of a hyperCR
Toda space, this Ricci-flat metric admits another compatible hypercomplex
structure with respect to which the symmetry \emph{is} triholomorphic,
and so we have nontrivial examples of selfdual spaces with two compatible
hypercomplex structures.
In summary, when $f=1$, the scalar-flat K\"ahler metric~\eqref{sfkm}
is hyperK\"ahler with a Killing vector and an additional
hypercomplex structure. The Einstein-Weyl space with geodesic
symmetry obtained from the conformal submersion $\xi$ is $S^3$ (with
symmetry given by a Hopf fibration).

\subsection{Selfdual Einstein metrics with negative scalar curvature}
\noindent\newline
These are the most interesting examples for us as they generalize diagram
(1.0). The article~\cite{Tod3} implies that, for any $a,b\in\R$,
the scalar-flat K\"ahler metric generated by the monopole
$a(1-\frac12zu_z)+\frac12bu_z$
on any Toda space is conformal to an
Einstein metric with scalar curvature $-3a$, via the conformal factor
$1/(az-b)^2$. For $a=0$ these are the hyperK\"ahler metrics discussed
in~\ref{HK} above, and for $a\neq0$ we can set $b=0$ by translating the $z$
coordinate. For our explicit solutions, this corresponds to adding a real
constant to $h$. Thus when $f=h$, i.e.,
$$\gmw=\frac{\tfrac12(h+\overline h)z+h\overline h}{(z+h)(z+\overline h)},$$
the scalar-flat K\"ahler metric~\eqref{sfkm} is conformal to an Einstein metric
with negative scalar curvature, given explicitly by
$$\frac1{z^2}\left[\frac{\tfrac12(h+\overline h)z+h\overline h}
{(z+h)(z+\overline h)}dz^2
+\bigl(\tfrac12(h+\overline h)z+h\overline h\bigr)g_{S^2}+
\frac{(z+h)(z+\overline h)}{\frac12(h+\overline h)z+h\overline h}
(dt+\theta)^2\right].$$
This has a conformal infinity at $z=0$ with conformal metric
$|h|^2g_{S^2}+\beta^2$, where $dt+\theta=\beta-\kappa z\,dz$ and so
$d\beta=\frac12(h+\overline h)\vol_{S^2}$.

Hence we see that we have found the selfdual Einstein metrics $M$ filling in
every Einstein-Weyl space admitting a geodesic symmetry.

We recover the Berger spheres by taking $h$ to be constant. The form of the
Einstein metrics we have found is easily related to the Pedersen family by
putting $h=1+im$ and setting $z=(1-\rho^2)/\rho^2$.

\section{Additional remarks}

We have shown that applying the LeBrun construction to an Einstein-Weyl space
with geodesic symmetry gives an Einstein metric with compatible hypercomplex
structure fibering over the general hyperCR Toda space.

The Einstein-Weyl spaces with geodesic symmetry are all hyperCR and this gives
an explanation for the hypercomplex structure coming with our Einstein
metrics. The twistor point of view gives a particularly quick way to see this:
a hyperCR structure on an Einstein-Weyl space $B$ corresponds to a holomorphic
map from its minitwistor space $\cS$ to $\CP1$. Composing this with the
projection from $PT\dual\cS$ to $\cS$, we see that the twistor space $Z$ of
$M$ has a holomorphic map to $\CP1$. Thus applying the LeBrun construction to
a hyperCR Einstein-Weyl space always gives a hypercomplex Einstein space.

There are many more hyperCR spaces than the spaces with geodesic symmetry
arising here, but it will be much harder to fill them in explicitly, since the
Einstein metric may no longer have a symmetry, so that it is harder to find
indirectly. It is perhaps easier to ask how other Einstein-Weyl spaces with
symmetry fill in. Indeed, it may be that there are other hyperCR spaces with
symmetry where the symmetry is \emph{not} geodesic, in which case we would
obtain selfdual Einstein metrics with a hypercomplex structure and a symmetry
which is not triholomorphic. However, even in this case, it is not clear how
to solve the Toda equation.

In summary, we note that we can generalize and augment diagram (1.0) as
follows:
\ifdiags
\begin{diagram}[width=1.5em,height=1.5em]
   &                      &\;M^4&       &\\
   &\ldTo^\xi[crab=4pt]\ruIntoB&&\rdTo^K&\\
B^3&                           &&       &\tilde B^3\\
   &\rdTo_K&&\ldTo_{\partial\rlap{$/\partial z$}}&\\
   &                      &S^2  &       &\\
\end{diagram}\fi
For the lower part of the diagram, we use the fact that the Killing field $K$
on $M^4$ descends to the geodesic symmetry of $B^3$, and $\partial/\partial z$
is a conformal submersion on the hyperCR Toda space $\tilde B^3$, since it is
shear-free. The surface over which $B^3$ and $\tilde B^3$ both fibre appears
to come with a natural spherical metric and so it would seem that it is the
geometry of $S^2$ which lies behind our constructions.

\end{document}